\documentclass{elsart3-1}

\usepackage[latin1]{inputenc}
\usepackage[T1]{fontenc}

\newcommand{\Z}{\mathbf{Z}}
\newcommand{\N}{\mathbf{N}}
\newcommand{\R}{\mathbf{R}}
\newcommand{\A}{\mathcal{A}}
\newcommand{\B}{\mathcal{B}}
\newcommand{\Bt}{\tilde{\mathcal{B}}}
\newcommand{\Ximin}{\Xi_\mathrm{min}}
\newcommand{\ra}{\rightarrow}
\newcommand{\omegat}{{\tilde{\omega}}}
\newcommand{\oubli}{\Theta}
\newcommand{\Xioubli}{\Xi_{\oubli}}
\newcommand{\D}{\mathcal D}

\usepackage{graphicx}

\usepackage{amssymb}

\usepackage[english,francais]{babel}

\newtheorem{theorem}{Theorem}[section]

\newtheorem{e-proposition}[theorem]{Proposition}

\newtheorem{e-definition}[theorem]{Definition\rm}

\newtheorem{theoreme}{Th\'eor\`eme}

\renewcommand{\theequation}{\arabic{equation}}
\def\og{\leavevmode\raise.3ex\hbox{$\scriptscriptstyle\langle\!\langle$~}}
\def\fg{\leavevmode\raise.3ex\hbox{~$\!\scriptscriptstyle\,\rangle\!\rangle$}}

\journal{the Acad\'emie des sciences}
\begin{document}
\centerline{}
\begin{frontmatter}


\selectlanguage{english}
\title{Combinatorics and topology of the Robinson tiling}


\selectlanguage{english}

\author[gaehler]{Franz G\"ahler},
\ead{gaehler@math.uni-bielefeld.de}
\author[julien]{Antoine Julien},
\ead{antoinej@uvic.ca}
\author[savinien]{Jean Savinien}.
\ead{savinien@univ-metz.fr}

\address[gaehler]{University of Bielefeld}
\address[julien]{University of Victoria--PIMS}
\address[savinien]{Universit\'e de Metz}



\begin{abstract}
\selectlanguage{english}
We study the space of all tilings which can be obtained
using the Robinson tiles (this is a two-dimensional subshift of finite type).
We prove that it has a unique minimal subshift, and describe it by means of a
substitution.
This description allows to compute its cohomology groups, and prove that it is
a model set.

\vskip 0.5\baselineskip

\selectlanguage{francais}
\noindent{\bf R\'esum\'e} \vskip 0.5\baselineskip \noindent
{\bf Combinatoire et topologie des pavages de Robinson.}
Nous \'etudions l'espace de tous les pavages qui peuvent s'obtenir
\`a partir des tuiles de Robinson (il s'agit d'un sous-d\'ecalage de type fini).
Cet espace contient un unique sous-espace minimal, que nous d\'ecrivons par le
biais d'une substitution.
En cons\'equence, il est possible de calculer les groupes de cohomologie associ\'es, et de montrer
qu'il s'agit d'un pavage de coupe et projection.
\end{abstract}
\end{frontmatter}

\selectlanguage{francais}
\section*{Version fran\c{c}aise abr\'eg\'ee}

C'est en 1971 que Robinson introduit l'ensemble de tuiles qui porte son nom. Un \og pavage de Robinson \fg{} est un pavage que l'on peut obtenir \`a partir des tuiles de la figure~\ref{fig:tiles-rob} (ainsi que leurs images par rotation et reflexion). Les pavages de Robinson doivent en outre respecter les r\`egles suivantes : les tuiles doivent se rencontrer face-\`a-face, et les fl\`eches doivent rencontrer des lignes ; par ailleurs, dans une colonne sur deux et une ligne sur deux, 
une tuile sur deux est de type (a) (voir fig.~\ref{fig:tiles-rob}), sans restriction \emph{a priori} sur son orientation.
Les tuiles de type (a) sont appel\'ees des \og carrefours \fg{}.

Formellement, un pavage est une d\'ecoration de $\Z^2$: \`a chaque \'el\'ement du r\'eseau correspond une tuile dans une orientation donn\'ee. Ainsi, un pavage est un \'el\'ement de $\A^{\Z^2}$, o\`u $\A$ est l'ensemble des tuiles de Robinson.
On note $\Xi$ l'ensemble des pavages de Robinson. C'est un sous-d\'ecalage de $\A^{\Z^2}$, c'est-\`a-dire un sous-ensemble ferm\'e (donc compact), et invariant sous l'action de $\Z^2$ par d\'ecalage (translation).
Un point important est que cet espace est non vide, et ne contient aucune p\'eriode (voir figure~\ref{fig:supertiles} pour un amas de taille~$7 \times 7$).

L'espace des pavages de Robinson n'est pas minimal, mais contient un unique sous-espace minimal, not\'e $\Ximin$.
Le th\'eor\`eme principal de cet article s'\'enonce ainsi.

\begin{theoreme}
Il existe une substitution (voir ci-dessous) $\tilde{\omega}$ telle que l'espace de pavages associ\'e \`a la substitution $\Xi_{\tilde{\omega}}$ est topologiquement conjugu\'e \`a $\Ximin$: il existe un hom\'eomorphisme entre ces deux espaces qui commute aux actions.
\end{theoreme}

Dans cet article, une substitution $\omega$ est une application qui associe 
a chaque tuile un carr\'e de $2 \times 2$ tuiles.
L'it\'eration d'une substitution produit une suite d'amas de taille croissante, et par un passage \`a la limite ad\'equat, produit un pavage $T$.
On appelle $\Xi_{\omega}$ le plus petit sous-d\'ecalage de $\B^{\Z^2}$ qui contient $T$.  Sous certaines conditions qui sont satisfaites ici, l'espace $\Xi_\omega$ ne d\'epend que de $\omega$, et pas de $T$.  Il est de plus minimal et sans p\'eriode.

La preuve du th\'eor\`eme ci-dessus se fait en deux parties.
D'abord, on exhibe une substitution $\omega$, d\'ecrite en figure~\ref{fig:subst1}.
Notons $\oubli$ l'application \og oubli \fg{} $\A \rightarrow \A \cup \{\mathrm{blanche}\}$, qui est l'identit\'e lorsque restreinte aux carrefours, et envoie toutes les autres tuiles sur la tuile blanche.  Alors cette application s'\'etend aux pavages, et d\'efinit une application $\Ximin \rightarrow \Xioubli$.
Alors $\Xi_\omega$ et $\Xioubli$ sont topologiquement conjugu\'es.
En d'autres termes, la substitution $\omega$ d\'ecrit enti\`erement l'emplacement des carrefours sur les pavages de Robinson.  Cela dit, certains pavages de Robinson ne sont pas enti\`erement d\'etermin\'es par la donn\'ee de tous leurs carrefours.
Ainsi, $\Xi_\omega$ n'est pas conjugu\'e \`a $\Ximin$.

La seconde \'etape est donc de d\'ecorer la substitution $\omega$ en une substitution \og augment\'ee \fg{}, not\'ee $\tilde{\omega}$, de sorte que les pavages obtenus par la nouvelle substitution d\'ecrivent l'emplacement des carrefours, mais aussi des autres tuiles qui composent les pavages de Robinson (les carrefours sont reli\'es par des lignes simples ou doubles).  Le processus de d\'ecoration est illustr\'e en partie par la figure~\ref{fig:single-line}, qui d\'ecrit comment encoder la position des lignes simples (ce sont les lignes fl\'ech\'ees des tuiles (b) et (c), et les lignes non fl\'ech\'ees des tuiles (b) et (c), voir figure~\ref{fig:tiles-rob}).

La conjugaison entre $\Ximin$ et $\Xi_{\tilde{\omega}}$ est construite explicitement par des d\'erivations locales. On donne une mani\`ere de re-coder un pavage de $\Ximin$ en terme des tuiles de $\tilde\omega$, et r\'eciproquement.

Ce r\'esultat a deux applications notables. Tout d'abord, une telle description permet de calculer des invariants topologiques pour l'espace $\Ximin$, en utilisant les m\'ethodes d'Anderson et Putnam~\cite{AP98}. Si on note $\Omega$ la suspension de $\Ximin$, on a:
 \[
   \check{\mathrm{H}}^2 (\Omega) = \Z[1/4] \oplus \big( \Z[1/2] \big)^{10} \oplus \Z^8 \oplus \Z / 4\Z; \quad
   \check{\mathrm{H}}^1 (\Omega) = \big( \Z[1/2] \big)^2 \oplus \Z ; \quad
   \check{\mathrm{H}}^0 (\Omega) = \Z.
 \]
Par ailleurs, comme l'espace (minimal) des pavages de Robinson est substitutif,
et chaque element contient un sous-ensemble p\'eriodique de carrefours,
on peut appliquer un th\'eor\`eme de Lee et Moody~\cite[th\'eor\`eme~3]{LM01}, qui 
implique que $\Ximin$ peut \^etre d\'ecrit par la m\'ethode de coupe et projection.
En particulier, le spectre de diffraction est purement ponctuel.

Avant de conclure, notons qu'une substitution plus simple, mais avec 
recouvrement de tuiles, et d\'ecrivant les m\^emes pavages de Robinson, 
a \'et\'e d\'ecouverte ind\'ependemment par Joan Taylor (communication priv\'ee).
Cette substitution permet aussi de calculer la cohomologie,
et donne les m\^emes r\'esultats.

\newpage

\selectlanguage{english}

\section{Introduction}

\begin{figure}[htp]
 \begin{center}
  \includegraphics[scale=0.75]{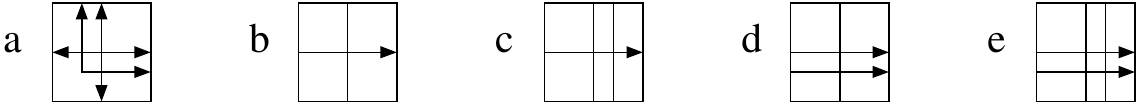}
 \end{center}
  \caption{The Robinson tiles. The tiles of the first type (a) are called \emph{crosses}.}
\label{fig:tiles-rob}
\end{figure}

In 1971, Robinson~\cite{Rob71} introduced his aperiodic set of tiles, in order
to build a two-dimensional subshift of finite type with no periodic orbit.
We also refer to the very comprehensive paper by Johnson and Madden~\cite{JM97}.
The Robinson tiling space is the set of all tilings which can be built from
the Robinson tiles as follows. Consider a set $\A$ of 28 symbols. These
symbols are the tiles of Figure~\ref{fig:tiles-rob}, as well as their images 
under rotations and reflections. Then, consider the subshift of finite type 
$\Xi \subset \A^{\Z^2}$ defined by:
\begin{enumerate}
 \item adjacency relations: two neighbouring tiles should meet in such a way
that arrowheads of a tile meet arrowtails of its neighbour;
 \item alternating cross rule: for all $x \in \Xi$, there is an element
$n \in \Z^2$, such that for all $i \in (2\Z) \oplus (2\Z)$, $x_{n+i}$ is
a \emph{cross} (see Figure~\ref{fig:tiles-rob}). There may also be crosses
at other positions.
\end{enumerate}

Then, it is possible to show that $\Xi$ is not empty, and none of its 
elements has periods under the action of $\Z^2$ by translation.
The key to this result is that any $x \in \Xi$ has a hierarchical structure:
for all $n\in \N$, define a $n$-supercross as shown in Figure~\ref{fig:supertiles}
(a $2$-supercross is given as an example, on the left).
Supercrosses are admissible under the adjacency rules. 
Therefore, by taking an appropriate union, one can
build an element of $\Xi$. Conversely, the matching rules force supercrosses
to appear in any admissible tiling. This gives a hierarchical structure which
allows to prove non-periodicity of any $x \in \Xi$.

An increasing union of $n$-supercrosses is called an \emph{infinite
order supertile}. Since there are several ways of including a $n$-supercross
in a $(n+1)$-supercross, there exist many different infinite order supertiles.
They need not cover the whole plane. Given $x \in \Xi$, we have the following
{\bf alternative}:
\begin{itemize}
 \item[(i)] either $x$ is made of only one infinite order supertile;
 \item[(ii)] or $x$ contains several (actually $2$ or $4$) infinite order
supertiles.
\end{itemize}

\begin{e-proposition}
The subshift $\Xi$ contains a unique, minimal subspace, called $\Ximin$.
Any element of $\Xi$ which follows alternative~(i) is in
$\Ximin$
\end{e-proposition}

To prove it, remark that any element of the tiling space contains
$n$-supercrosses for all $n$. Therefore, elements which follow alternative~(i) are accumulation
points of any orbit of $\Xi$.

\section{A substitution}

We want to define a substitution map (see for example~\cite{AP98}) which
describes the Robinson minimal tiling space.
In our context, a substitution is a set of tiles $\B$, and a map $\omega:
\B \ra \B^{\{0,1\}\times \{0,1\}}$ which associates to every tile  a $2 \times 2$
patch of tiles.
Then the substitution tiling space $\Xi_\omega$ is the set of all elements $x \in \B^{\Z^2}$
such that any patch $x_{[n,n+k)\times[m,m+l)}$ of any size $k \times l$ appears
(up to translation) as a subpatch of $\omega^N(t)$ for some integer $N$ and some
tile $t$.
It is a closed, shift-invariant subset.

\begin{figure}[htp]
 \begin{center}
  \includegraphics[scale=0.75]{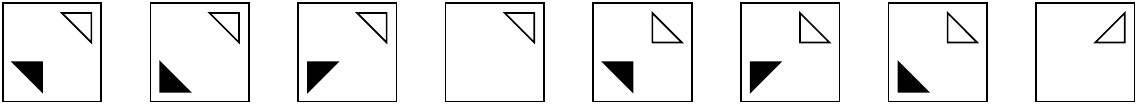}
 \end{center}
\caption{The tiles of the substitution.}
\label{fig:tuiles-subst1}
\end{figure}

\begin{figure}[htp]
 \begin{center}
  \includegraphics[scale=0.75]{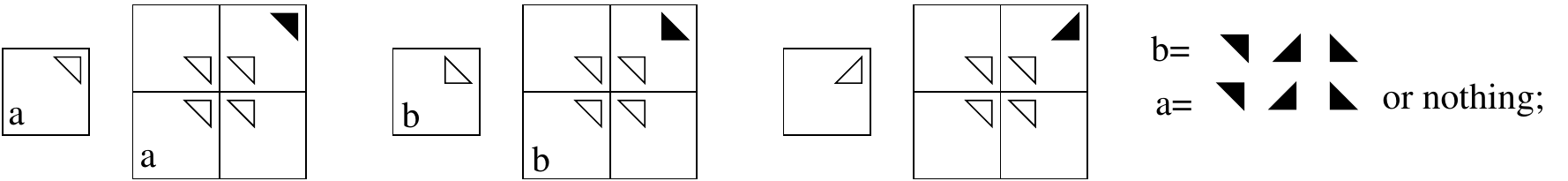}
 \end{center}
\caption{The substitution.}
\label{fig:subst1}
\end{figure}

Consider the tiles given in Figure~\ref{fig:tuiles-subst1}, and the associated
substitution given by Figure~\ref{fig:subst1}.
One can check that this substitution is \emph{primitive}: for
any two tiles $a$ and $b$, the tile $b$ is contained in the $n$-th substitution
of $a$, for $n$ big enough.
It is known that tiling spaces associated with primitive substitutions
are \emph{minimal}.

\begin{theorem}\label{thm:factor}
The tiling space $\Xi_\omega$ associated with the substitution above is a factor
of $\Ximin$, that is, there is an onto map $\phi: \Ximin \ra \Xi_\omega$, which
commutes with the shift.
\end{theorem}

We describe this factor more precisely in the next section. Then, we will prove
the theorem.
Let us mention without proof that this factor map is almost everywhere one-to-one.
Moreover, we will explain below, by ``decorating'' the tiles and the 
substitution, how one can build a substitution $\tilde{\omega}$, 
whose tiling space is topologically conjugate to $\Ximin$.

\section{Local derivations}

A map $\phi: \A^{\Z^2} \ra \B^{\Z^2}$ is called a \emph{local derivation}
(see~\cite{BSJ}) if for any $x$, and any position $(i,j)$, the tile
$\phi(x)_{i,j}$ only depends on the configuration of $x$ \emph{around} $(i,j)$.
That is, it only depends on the $x_{k,l}$, for $\|(k-i,l-j)\| < C$, where $C$
only depends on $\phi$.
A local derivation is automatically continuous, and commutes with the shift.

Consider two minimal subshifts $\Xi_\A \subset \A^{\Z^2}$ and $\Xi_\B \subset \B^{\Z^2}$,
and a local derivation $\phi$.
If $x \in \Xi_\A$ and $\phi(x) \in \Xi_\B$, then
by minimality, $\phi$ maps $\Xi_\A$ onto $\Xi_\B$, and therefore,
it is a factor map.
If there are two local derivations $\Xi_\A \ra \Xi_\B$ and $\Xi_\B \ra \Xi_\A$
which are inverse of each other, then the two subshifts are topologically
conjugate.
In order to prove Theorem~\ref{thm:factor}, we need to define a local derivation
$\Ximin \ra \B^{\Z^2}$, and then prove that at least one point of $\Ximin$
is mapped to $\Xi_\omega$.

\begin{figure}[htp]
 \begin{center}
  \includegraphics[scale=0.75]{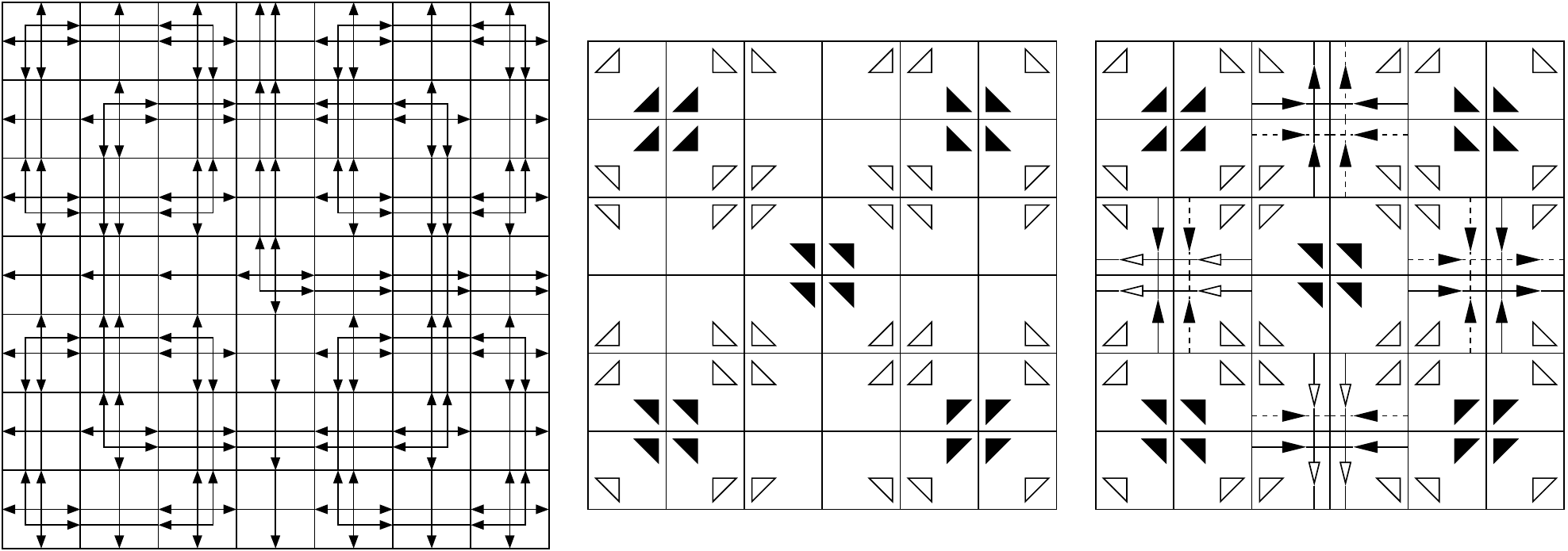}
 \end{center}
\caption{On the left, a $2$-super-cross: it is made of four 1-super-crosses, and
one cross in the middle. The middle cross is connected to lines which extend
to the boundary of the super-cross: the supercross behaves like a cross of
bigger size.
In the middle, the image of this super-cross by the local derivation $\phi$
is shown.
On the right, a fully decorated version of this super-cross (see Section~\ref{sec:decoration}) is given.}
\label{fig:supertiles}
\end{figure}

This local derivation is a composition of two maps $\phi = \phi_2 \circ \phi_1$.
The map $\phi_1$ is a ``forgetful'' map, which only remembers the position and
orientation of the crosses.
The map $\phi_2$ breaks the tiles and recomposes them.
An example of a patch together with its image by $\phi$ is given by
Figure~\ref{fig:supertiles}.
More precisely, $\phi_1$ is defined as follows (note that whether a cross is
subject or not to the alternating cross rule is a local information):
\begin{itemize}
 \item The image of a cross which is subject to the alternating cross rule is
a tile decorated with empty triangles which remember the orientation of the cross;
 \item The image of any other cross is a tile decorated with solid black triangles, which
remember the orientation of the cross;
 \item The image of any other tile is a blank (empty) tile.
\end{itemize}

Remark that in a super-cross, every row and every column contains
at least one cross.
Therefore, all the tiles in a super-cross are entirely determined by the 
positions of the crosses.
For this reason, the ``forgetful'' derivation described above is one-to-one
when restricted to the set of Robinson tilings which follow alternative~(i).

The map $\phi_2$ is defined by cutting all tiles obtained above in four,
and recomposing them, so that the new tiles are now made of four pieces of four
different previous tiles.
The set of tiles which can be obtained is exactly $\B$, on which the substitution
is defined.
This derivation is of course invertible, by cutting tiles again and recomposing them (see Figure~\ref{fig:supertiles}).
We call $\phi$ the composition of $\phi_1$ and $\phi_2$.

In order to prove that $\phi$ maps $\Ximin$ to $\Xi_\omega$, we need to iterate
the substitution $\omega$ on a tile. Remark that for any tile $t$, $\omega(t)$ contains
the image of a cross under the local derivation.
By iteration, one proves that $\omega^n(t)$ contains a $(n-1)$-super-cross
($0$-super-crosses being simply crosses). Therefore, the image under the local
derivation of any Robinson tiling made of a single infinite order super-cross is
contained in $\Xi_\omega$. We conclude, using minimality, that $\Xi_\omega$ is a factor
of~$\Ximin$.

\section{Decorating the substitution}
\label{sec:decoration}

The idea is now to decorate tiles of $\B$ to get a set of tiles $\Bt$, and a substitution
$\omegat$ on it.
If we call $\psi$ the map $\Bt \ra \B$ which forgets the decorations, we require
that $\psi \circ \omegat = \omega \circ \psi$ (so that $\omega$ is induced by
$\omegat$ on undecorated tiles).
Then, $\psi$ induces a local derivation $\Xi_{\omegat} \ra \Xi_\omega$. 
The augmented substitution $\omegat$ should be defined in such a way that the two spaces $\Ximin$ and
$\Xi_\omegat$---of which $\Xi_\omega$ is a factor---are conjugate.

Note that $\Xi_\omega$ and $\Ximin$ cannot be conjugate to each other:
tilings containing more than one infinite order supertile in $\Ximin$
are mapped to tilings in $\Xi_\omega$ which may have several pre-images.
The reason is that in such a tiling, there may be one row or one column in which
there is no cross. The map $\phi$ forgets everything about such a line.

In order to get a space which is conjugate to $\Ximin$, we need to keep track
in $\omegat$ of the lines which connect crosses (simple
lines, double lines, together with arrows orienting them).

\paragraph*{Keeping track of simple lines.}
We describe how to produce decorations to keep track of simple
lines (arrowed lines in tiles of type~(b) and~(c), and non-arrowed lines in tiles
of type~(b) and~(d), see Fig.~\ref{fig:tiles-rob}).
\begin{figure}[htp]
 \begin{center}
  \includegraphics[scale=0.75]{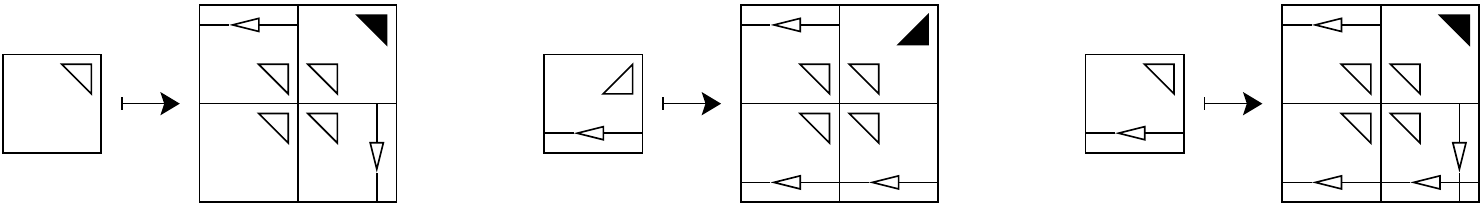}
 \end{center}
\caption{Single lines are coded by white-arrowed lines.
In the first substitution, arrows were added near two edges. If such a patch
occurs in a tiling $y \in \Xi_\omega$, then any $x \in \Ximin$ which maps
to $y$ will have a cross at the position of the black triangle. So, the 
white-arrowed lines in $y$ correspond to single lines in $x$.
This creates new tiles. In the middle, the substitution
of the lower-right tile of the first patch is shown,
and on the right, the substitution of the lower-left tile of the second
patch.}
\label{fig:single-line}
\end{figure}
Simple lines connect crosses which are ``back-to-back''.
Therefore, we add decorations to the tiles in $\Bt$ in order to remember this
fact (see Figure~\ref{fig:single-line}).
Next, we decide how to substitute a decorated tile. This creates new 
decorated tiles which are not yet in $\Bt$. 
Then, we decide how these decorated tiles substitute.
We iterate this process until no new tile is created. Finally, we discard any tile which
doesn't appear any more in the eventual range of our decorated substitution.

In order to have a fully decorated substitution, we also need to keep track of
double lines. See Figure~\ref{fig:supertiles} (patch on the right) for an example
of a fully decorated patch. 
The decorated substitution $\tilde\omega$ requires $208$ different tiles, 
up to translation.

\begin{theorem}
 Call $\tilde{\omega}$ the substitution defined on decorated tiles.
Then the subshifts $\Ximin$ and $\Xi_{\tilde{\omega}}$ are topologically
conjugate.
\end{theorem}

To prove this, we have to find rules to associate a tile of the Robinson tiling
to any $2 \times 2$ patch of decorated tiles. This gives a local derivation which
induces a factor map $\Xi_{\tilde{\omega}} \rightarrow \Ximin$.
In the converse direction, we have to find a derivation from $2 \times 2$ Robinson patches to decorated tiles.
The two derivations must be inverses of each other, in the sense that the 
image of any admissible $3 \times 3$ patch of a Robinson tiling under the
composition of the two derivations is the central tile.

\section{Applications}

Since we now have a substitution $\tilde{\omega}$ on $\tilde{\B}$, it is possible
to use the methods developed in~\cite{AP98} to compute
the cohomology of the continuous hull of $\Ximin$.
By continuous hull (or tiling space), we mean the suspension of the action of $\Z^2$ on $\Ximin$:
if $\sigma^{(n,m)}(x)$ is the image of $x$ by the shift, define:
\[
 \Omega = \big( \Ximin \times \R^2 \big) / \big\{ (x,(t_1, t_2)) \sim (\sigma^{(n,m)}(x), (t_1 - n, t_2 - m)) \big\}.
\]
Given a substitution, one can associate a finite CW-complex $\Gamma$, with $2$-cells being
(the interior of) the tiles. Then, the tiling space is homeomorphic to the inverse
limit of $\Gamma$ under a map induced by the substitution.
Using a computer program, it was possible to determine this complex,
(adjacencies, etc.), compute its \v{C}ech cohomology, and the map induced by
the substitution on cohomology.

\begin{theorem}
The cohomology groups of the Robinson minimal tiling space are:
 \[
   \mathrm{H}^2 (\Omega) = \Z[1/4] \oplus \big( \Z[1/2] \big)^{10} \oplus \Z^8 \oplus \Z / 4\Z; \quad
   \mathrm{H}^1 (\Omega) = \big( \Z[1/2] \big)^2 \oplus \Z ; \quad
   \mathrm{H}^0 (\Omega) = \Z.
 \]
\end{theorem}

Another consequence of the fact that the Robinson tiling can be described by
a (lattice) substitution is that it is a model set. As every element of the 
hull contains a lattice periodic subset of crosses, this follows from
Theorem~3 of~\cite{LM01}. As a further consequence, it also implies
that every tiling in the hull has pure point diffraction spectrum.

Before concluding, let us mention that a somewhat simpler substitution, 
although with overlapping tiles, has been found independently by Joan Taylor
(private communication).
Also that substitution can be used to compute the cohomology, and gives
the same results.

\clearpage

\appendix

\section{Complements to the article}

This appendix presents arguments (and pictures) for several results which were outlined in the article.

\subsection{Non-minimality of the Robinson space}

\begin{figure}[htp]
 \begin{center}
  \includegraphics[scale=0.75]{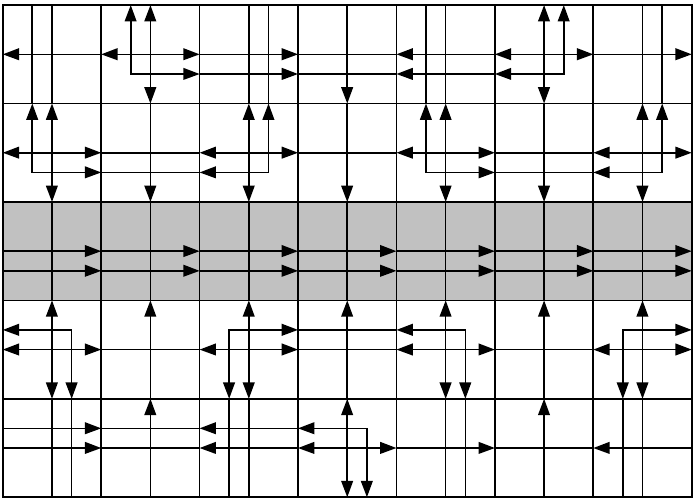}
 \end{center}
\caption{The grey line extends to infinity on both sides. On top and bottom of
it are two infinite order supertiles. One can shear the top and bottom infinite
order supertile independently by any multiple of $2$ and still get a Robinson tiling.
However, such a tiling is not in $\Ximin$.}
\label{fig:not-minimal}
\end{figure}

Non-minimality of the set of all admissible Robinson tilings is illustrated by Figure~\ref{fig:not-minimal}.
It shows a Robinson tiling made of two infinite order supertiles, separated by an infinitely long line of tiles.
The two infinite order supertiles can be ``sheared'': if one translates the bottom infinite-order supertile to the right by any multiple of two, the result is still admissible (and satisfies the alternating cross rule). However, such a sheared tiling is not in $\Ximin$: the patch of Figure~\ref{fig:not-minimal} cannot appear in a super-cross.

\subsection{The factor map $\Ximin \ra \Xi_\omega$ is non-trivial}

\begin{figure}[htp]
 \begin{center}
  \includegraphics[scale=0.75]{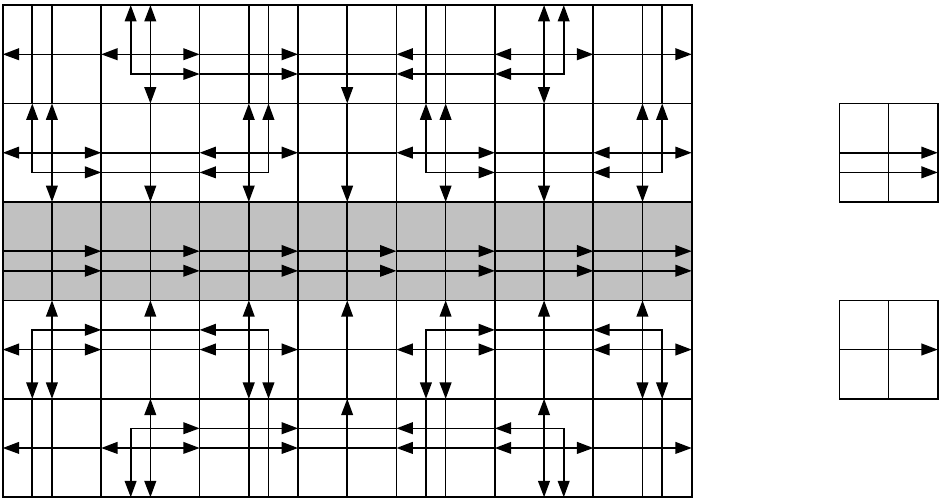}
 \end{center}
\caption{The grey line extends to infinity on both sides. On top and bottom of
it are two infinite order supertiles. If one changes the tiles on the grey line
by any of the two tiles on the right (or their reflections), one
gets different Robinson tilings, but their images in $\Xi_\omega$ are 
nevertheless the same.}
\label{fig:not-one-to-one}
\end{figure}

We present several Robinson tilings which have the same image in $\Xi_\omega$.
Such tilings are made of two infinite order supertiles, separated by an infinitely long line of tiles (of thickness $1$ tile).
This is described by Figure~\ref{fig:not-one-to-one}.

In this picture, the grey row of tiles contains no cross. Therefore, its image in $\Xi_\omega$ is a row of blank tiles.
If one changed the tiles which compose the grey row, the image would still be the same.
Actually, one can check that there are six such Robinson tilings which have the same image in $\Xi_\omega$.

\subsection{Decorated substitution}

\begin{figure}[htp]
 \begin{center}
  \includegraphics[scale=0.75]{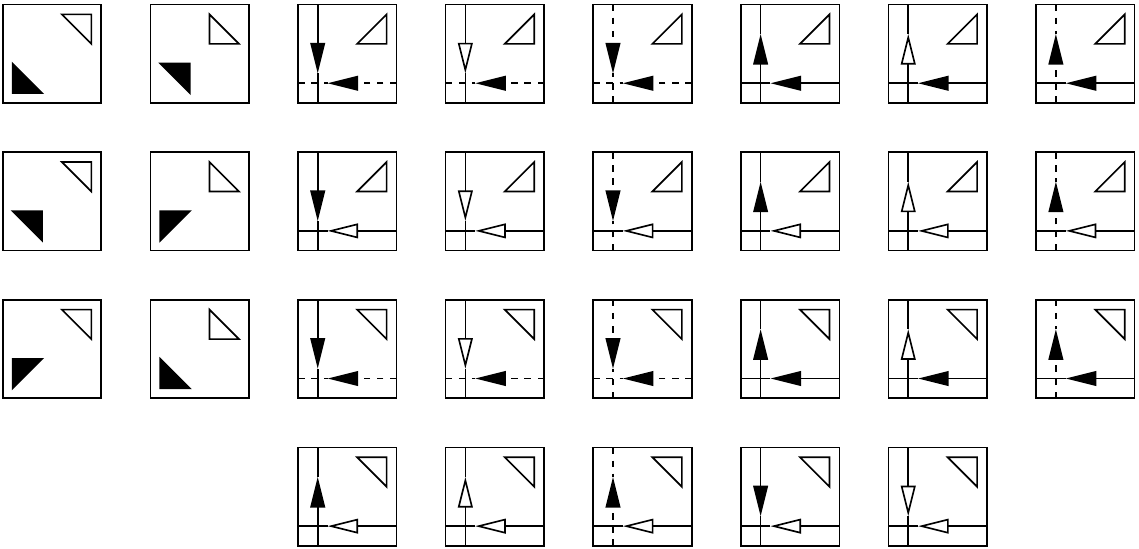}
 \end{center}
\caption{Decorated tiles of $\omegat$.}
\label{fig:decorated-tiles}
\end{figure}

\begin{figure}[htp]
 \begin{center}
  \includegraphics[scale=0.75]{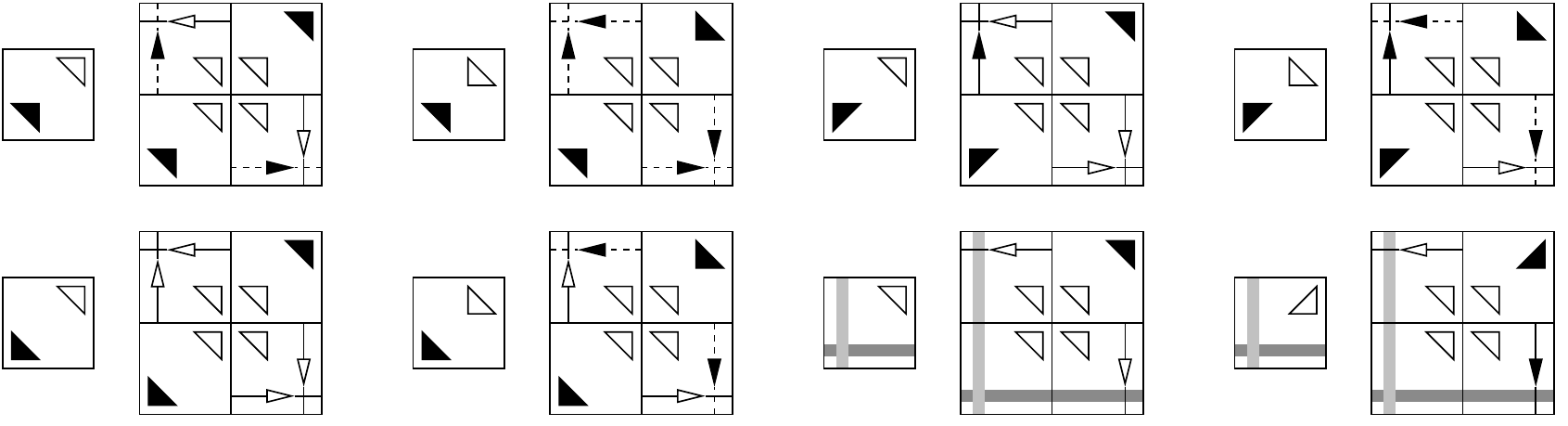}
 \end{center}
\caption{The decorated substitution $\omegat$.
The two types of grey lines can be (possibly different) arrows of any type,
\emph{a priori} (black on plain line, black on dotted line or white, in any
direction).
The vertical (resp.\ horizontal) arrows in the substitution of a tile are
of the same kind as the vertical (resp.\ horizontal) arrows in the tile.
Note also that \emph{a posteriori}, not all combinations of arrows are
acceptable. See Figure~\ref{fig:decorated-tiles} for a list of all legal
tiles.}
\label{fig:decorated-subst}
\end{figure}

We present here the decorated substitution. It is defined on the set of tiles given by Figure~\ref{fig:decorated-tiles} (as well
as their images under rotation and reflection).
The substitution itself is described by Figure~\ref{fig:decorated-subst}.

A direct check shows that forgetting the decorations gives back the undecorated substitution $\omega$. It was checked, using a computer, that the substitution $\omegat$ is primitive.
As we already described, the augmented substitution appears to be a superposition of $\omega$, and of a one-dimensional substitution on arrows: the image of an arrow is two copies of itself (symbolically: $a \mapsto aa$).

\subsection{Local derivations}
The local derivations from $\Ximin$ to $\Xi_{\tilde\omega}$ and back 
can now easily be derived.
Call $\D_1$ the local derivation $\Ximin \ra \Xi_{\tilde\omega}$, and $\D_2$ 
the local derivation in the other direction. We describe these now.
As a preliminary remark, we note that tiles of the Robinson tiling which are 
not crosses are defined by two features: vertical lines, and horizontal lines.
Similarly, decorated tiles in $\Bt$ can be seen as the superposition of one, two or three features: a tile in $\B$, and possibly a vertical arrow, and possibly a horizontal arrow.

\begin{figure}[htp]
 \begin{center}
  \includegraphics[scale=0.75]{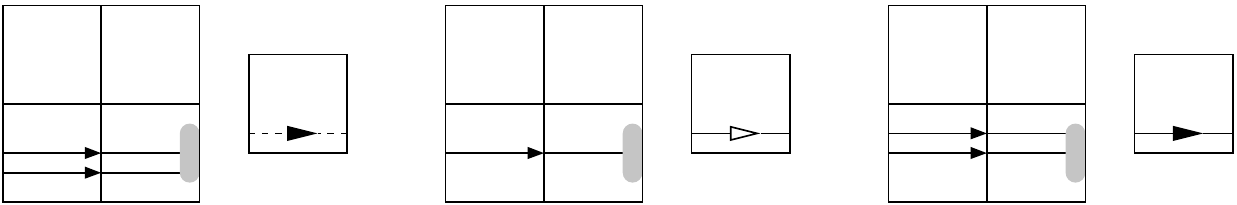}
 \end{center}
\caption{The local derivation adding a decoration to the tiles of $\B$. If a $2\times 2$ patch has an indicated feature (regardless of what the other tiles or vertical lines are), then the image tile has the indicated decoration. The grey ends of the lines can be arrowed, but don't need to: the image by the derivation is the same.
Note that the rules for adding vertical decorations on the tiles can be deduced by rotating these pictures.}
\label{fig:deriv1}
\end{figure}

Let us first describe $\D_1$: it is a map from $\A^{\{0,1\}\times \{0,1\}}$ to $\Bt$. We already described a map to $\B$, it is enough to describe how to add decorating features. This is done in Figure~\ref{fig:deriv1}.

\begin{figure}[htp]
 \begin{center}
  \includegraphics[scale=0.75]{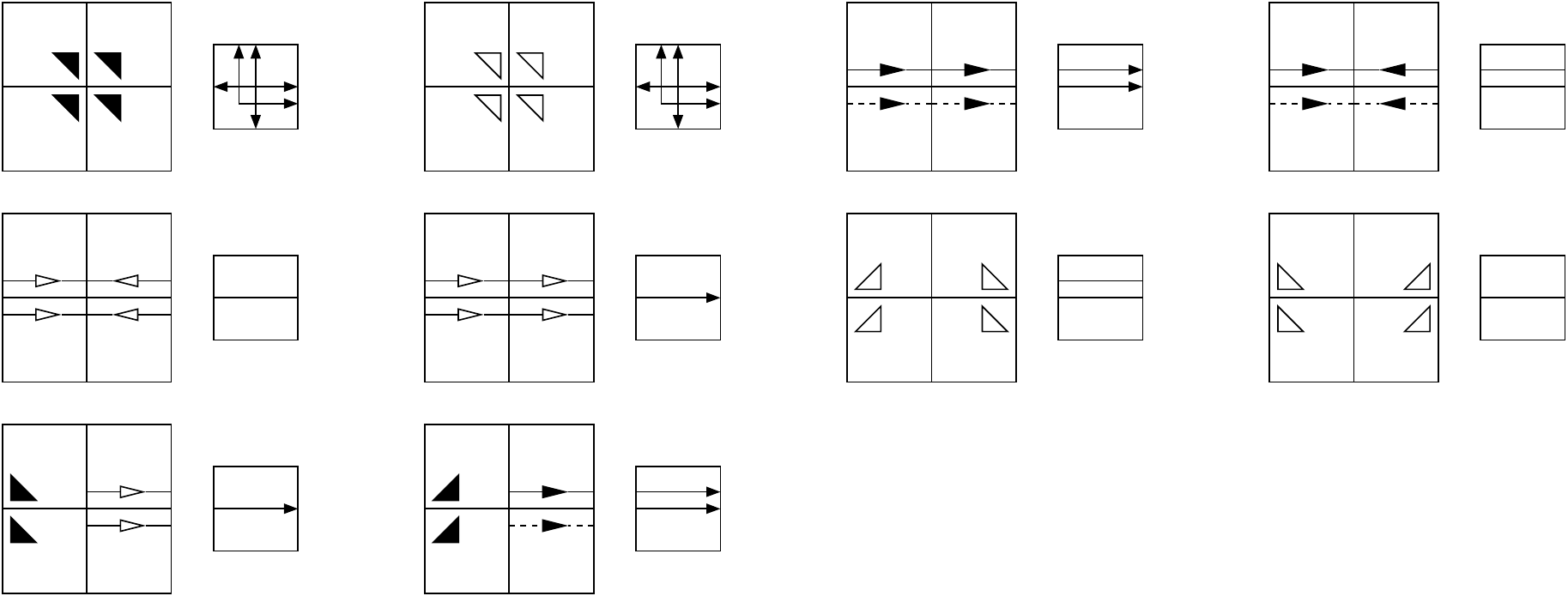}
 \end{center}
\caption{The local derivation $\D_2$. If a $2 \times 2$ patch has one of the indicated features, then its image has the indicated horizontal lines.
As in Figure~\ref{fig:deriv1}, the rules for the vertical lines can be obtained by rotating the pictures.}
\label{fig:deriv2}
\end{figure}

The converse derivation $\D_2$ is defined similarly by a map from $\Bt^{\{0,1\}\times \{0,1\}}$ to $\A$. It is described in Figure~\ref{fig:deriv2}.

\begin{figure}[htp]
 \begin{center}
  \includegraphics[scale=0.75]{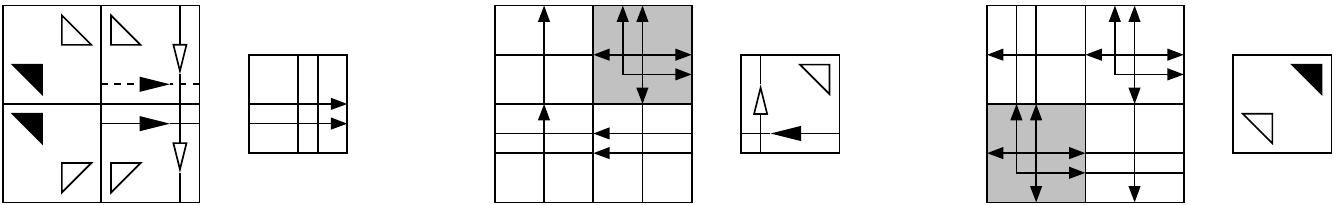}
 \end{center}
\caption{Three examples of local derivations using the rules above. The tiles subject to the alternating cross rule were previously given a grey background (this change is a local derivation itself).}
\label{fig:local-deriv-ex}
\end{figure}

\begin{figure}[htp]
 \begin{center}
  \includegraphics[scale=0.75]{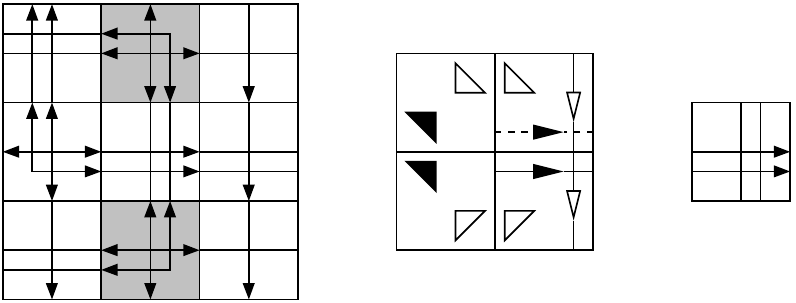}
 \end{center}
\caption{The local derivations applied back and forth. Note that the final tile is the same as the central tile of the starting patch.}
\label{fig:local-deriv-ex-itere}
\end{figure}

Figure~\ref{fig:local-deriv-ex} shows three examples of local derivations using the rules above. Figure~\ref{fig:local-deriv-ex-itere} shows a local derivation on a $3 \times 3$ patch, composed with the local derivation backwards. The result is a single tile.

More generally, if the composition of these two derivations is applied on \emph{any} $3 \times 3$ patch, the image is the centre tile of the original patch. This fact was checked on a computer.
This proves that the factor maps induced by these derivations $\Ximin \ra \Xi_\omegat$ and $\Xi_\omegat \ra \Ximin$ are inverse of each other. Therefore, the spaces are topologically conjugate.


\end{document}